\documentclass[12pt]{article}
\usepackage{amsfonts,amsmath,amsxtra}
\usepackage{amssymb, amscd}
\def\hybrid{\topmargin 0pt      \oddsidemargin 0pt
        \headheight 0pt \headsep 0pt
        \textwidth 16.5cm
        \textheight 23cm
        \marginparwidth 0.0in
        \parskip 5pt plus 1pt   \jot = 1.5ex}
\catcode`\@=11
\def\marginnote#1{}
\newcount\hour
\newcount\minute
\newtoks\amorpm
\hour=\time\divide\hour by60
\minute=\time{\multiply\hour by60 \global\advance\minute by-\hour}
\edef\standardtime{{\ifnum\hour<12 \global\amorpm={am}%
        \else\global\amorpm={pm}\advance\hour by-12 \fi
        \ifnum\hour=0 \hour=12 \fi
      \number\hour:\ifnum\minute<10 0\fi\number\minute\the\amorpm}}
\edef\militarytime{\number\hour:\ifnum\minute<10 0\fi\number\minute}

\def\draftlabel#1{{\@bsphack\if@filesw {\let\thepage\relax
   \xdef\@gtempa{\write\@auxout{\string
      \newlabel{#1}{{\@currentlabel}{\thepage}}}}}\@gtempa
   \if@nobreak \ifvmode\nobreak\fi\fi\fi\@esphack}
        \gdef\@eqnlabel{#1}}
\def\@eqnlabel{}
\def\@vacuum{}
\def\draftmarginnote#1{\marginpar{\raggedright\scriptsize\tt#1}}

\def\draft{\oddsidemargin -0.1truein
        \def\@oddfoot{\sl preliminary draft \hfil
        \rm\thepage\hfil\sl\today\quad\militarytime}
        \let\@evenfoot\@oddfoot \overfullrule 3pt
        \let\label=\draftlabel
        \let\marginnote=\draftmarginnote
\def\@eqnnum{{\rm (\theequation)}
\rlap{\kern\marginparsep\tt\@eqnlabel}%
\global\let\@eqnlabel\@vacuum}  }

\newcounter{app}

\def\app{\setcounter{equation}{0}
\def\theequation{\Alph{app}.\arabic{equation}}\par
   \addvspace{10ex}
   \@afterindentfalse
  \secdef\@app\@dapp}
\newcommand\@app{\@startsection {app}{1}{-0.3ex}%
                             {-3.5ex \@plus -1ex \@minus -.2ex}%
                                   {2.3ex \@plus.2ex}%
                                   {\normalfont\Large\bf}}

\def\@dapp#1{%
{\parindent \z@ \raggedright  \bf #1}\par\nobreak}
\def\l@app#1#2{\ifnum \c@tocdepth >\z@
    \addpenalty\@secpenalty
    \addvspace{1.0em \@plus\p@}%
    \setlength\@tempdima{1.5em}%
    \begingroup
      \parindent \z@ \rightskip \@pnumwidth
      \parfillskip -\@pnumwidth
      \leavevmode \bfseries
      \advance\leftskip\@tempdima
      \hskip -\leftskip
      #1\nobreak\hfil \nobreak\hb@xt@\@pnumwidth{\hss #2}\par
    \endgroup\fi}
\newcounter{sapp}[app]

\def\sapp{\def\theequation{\Alph{app}.\arabic{equation}}\par
   \@afterindentfalse
  \secdef\@sapp\@dsapp}
\newcommand\@sapp{\@startsection{sapp}{2}{\z@}%
                           {-3.25ex\@plus -1ex \@minus -.2ex}%
                           {1.5ex \@plus .2ex}%
                              {\normalfont\large\bfseries}}

\def\@dsapp#1{%
{\parindent \z@ \raggedright  \bf #1}\par\nobreak}

\newcommand{\l@sapp}{\@dottedtocline{2}{1.4em}{2.5em}}

\newcommand{\RR}{{\mathbb{R}}}
\newcommand{\CC}{{\mathbb{C}}}
\newcommand{\NN}{{\mathbb{N}}}       

\newcommand{\ZZ}{{\mathbb{Z}}}
\newcommand{\TT}{{\mathbb{T}}}

\newfont{\Bbbb}{msbm7 scaled 1\@ptsize00}
\newcommand{\zs}{\raise-1pt\hbox{$\mbox{\Bbbb Z}$}}

\font\sevenmsa=msam6 
\newfam\msafam
\textfont\msafam=\sevenmsa
\def\hexnumber@#1{\ifnum#1<10 \number#1\else
\ifnum#1=10 A\else\ifnum#1=11 B\else\ifnum#1=12 C\else
\ifnum#1=13 D\else\ifnum#1=14 E\else\ifnum#1=15 F\fi\fi\fi\fi\fi\fi\fi}
\def\msa@{\hexnumber@\msafam}
\def\llcorner{\delimiter"4\msa@78\msa@78 }
\def\lrcorner{\delimiter"5\msa@79\msa@79 }
\mathchardef\blacktriangleright="3\msa@49
\mathchardef\blacktriangleleft="3\msa@4A
\font\tenmsb=msbm10 scaled 1\@ptsize00
\newfam\msbfam
\textfont\msbfam=\tenmsb
\scriptfont\msbfam=\tenmsb

\newdimen\linethick  \linethick=0.4pt
\newdimen\hboxitspace    \hboxitspace=5pt
\newdimen\vboxitspace    \vboxitspace=5pt

\def\fr#1{%
\be\new
\vcenter{
\hrule height\linethick
           \hbox{\vrule width\linethick
                 \kern\hboxitspace
                 \vbox{\kern\vboxitspace
                       \hbox{$\begin{array}{c}\displaystyle#1
          \end{array}$}%
                       \kern\vboxitspace}%
                 \kern\hboxitspace
                 \vrule width\linethick}%
           \hrule height\linethick}%
\ee}

\newdimen\Squaresize \Squaresize=14pt
\newdimen\Thickness \Thickness=0.5pt

\def\Square#1{\hbox{\vrule width \Thickness
   \vbox to \Squaresize{\hrule height \Thickness\vss
      \hbox to \Squaresize{\hss#1\hss}
   \vss\hrule height\Thickness}
\unskip\vrule width \Thickness}
\kern-\Thickness}

\def\Vsquare#1{\vbox{\Square{$#1$}}\kern-\Thickness}

\def\numberbysection{\@addtoreset{equation}{section}
        \def\theequation{\thesection.\arabic{equation}}}
\numberbysection

\renewcommand{\theequation}{\thesection.\arabic{equation}}
\def\titlepage{\@restonecolfalse\if@twocolumn\@restonecoltrue\onecolumn
     \else \newpage \fi \thispagestyle{empty}\c@page\z@
        \def\thefootnote{\fnsymbol{footnote}} }

\def\endtitlepage{\if@restonecol\twocolumn \else  \fi
        \def\thefootnote{\arabic{footnote}}
        \setcounter{footnote}{0}}  
\relax

\hybrid
\makeatletter
\newdimen\normalarrayskip            
\newdimen\minarrayskip               
\normalarrayskip\baselineskip
\minarrayskip\jot
\newif\ifold             \oldtrue            \def\new{\oldfalse}
\def\arraymode{\ifold\relax\else\displaystyle\fi}
\def\eqnumphantom{\phantom{(\theequation)}} 
\def\@arrayskip{\ifold\baselineskip\z@\lineskip\z@
     \else
     \baselineskip\minarrayskip\lineskip1\baselineskip\fi}

\def\@arrayclassz{\ifcase \@lastchclass \@acolampacol \or
\@ampacol \or \or \or \@addamp \or
   \@acolampacol \or \@firstampfalse \@acol \fi
\edef\@preamble{\@preamble
  \ifcase \@chnum
     \hfil$\relax\arraymode\@sharp$\hfil
     \or $\relax\arraymode\@sharp$\hfil
     \or \hfil$\relax\arraymode\@sharp$\fi}}

\def\@array[#1]#2{\setbox\@arstrutbox=\hbox{\vrule
     height\arraystretch \ht\strutbox
     depth\arraystretch \dp\strutbox
width\z@}\@mkpream{#2}\edef\@preamble{\halign \noexpand\@halignto
\bgroup \tabskip\z@ \@arstrut \@preamble \tabskip\z@ \cr}%
\let\@startpbox\@@startpbox \let\@endpbox\@@endpbox
  \if #1t\vtop \else \if#1b\vbox \else \vcenter \fi\fi
  \bgroup \let\par\relax
  \let\@sharp##\let\protect\relax
  \@arrayskip\@preamble}
\def\eqnarray{\stepcounter{equation}%
              \let\@currentlabel=\theequation
              \global\@eqnswtrue
              \global\@eqcnt\z@
              \tabskip\@centering              
              \let\\=\@eqncr
              $$%
            \halign to \displaywidth  \bgroup
             \eqnumphantom \@eqnsel
      \hskip\@centering                               
    $\displaystyle  \tabskip\z@ {##}$%
    &\global\@eqcnt\@ne \hskip 2\arraycolsep
         $ \displaystyle  \arraymode{##}$\hfil
    &\global\@eqcnt\tw@ \hskip 2\arraycolsep
         $\displaystyle\tabskip\z@{##}$\hfil
         \tabskip\@centering
    &{##}\tabskip\z@\cr}
\makeatother
\newtheorem{te}{Theorem}[section]
\newtheorem{de}{Definition}[section]
\newtheorem{prop}{Proposition}[section]           

\newtheorem{rem}{Remark}[section]
\newcommand{\beq}[1]{\begin{equation}\label{#1}}
\newcommand\eeq{\end{equation}}
\newcommand\bqa{\begin{eqnarray}}
\newcommand\eqa{\end{eqnarray}}

\def\be{\begin{eqnarray}\new\begin{array}{cc}}
\def\ee{\end{array}\end{eqnarray}}

\def\beq{\begin{equation}}
\def\eeq{\end{equation}}
\def\bse{\begin{subequations}}                
\def\ese{\end{subequations}}
\def\bp{\begin{pmatrix}}
\def\ep{\end{pmatrix}}

\def\h{\hbar}
\def\i{\imath}

\renewcommand{\theequation}{\thesection.\arabic{equation}}

\def\d{\partial}
\def\stack#1#2{\raise0.7pt\hbox{$\mathrel{\mathop{#2}\limits^{#1}}$}}
\def\tr{\triangleright}
\def\tl{\triangleleft}
\def\sem{\mathsurround=0pt \raise1pt
\hbox{$\scriptscriptstyle>\!\!$}\:\!\!\tl}
\def\mes{\mathsurround=0pt \tr\!\:\!\raise0.8pt
\hbox{$\scriptscriptstyle\!\!<$}\,}
\def\]{\mathsurround=0pt ]\raise-2pt\hbox{$_\ast$}}

\def\al{\alpha}

\def\<{\langle}
\def\>{\rangle}

\def\frak{\mathfrak}

\def\ts#1#2{{\textstyle\frac{#1}{#2}}}

\def\we{\raise-1pt\hbox{$\,\stackrel{\wedge}{,}\,$}}

\def\b{\beta}

0
\begin{document}

\thispagestyle{empty}
\begin{center}

\phantom.
\bigskip\bigskip\bigskip\bigskip
{\hfill{\normalsize ITEP-TH-31/04}\\
\hfill{\normalsize HMI-04-06}\\
\hfill{\normalsize TCD-MATH-05-01}\\
[10mm]\Large\bf
On a Class of Representations of Quantum Groups}
\vspace{0.5cm}

\bigskip\bigskip
{\large A. Gerasimov}
\\ \bigskip
{\it Institute for Theoretical \& Experimental Physics, 117259,
Moscow, Russia}\\ {\it
 Department of Pure and Applied Mathematics, Trinity
College, Dublin 2, Ireland } \\ {\it Hamilton
Mathematics Institute, TCD, Dublin 2, Ireland}\\
\bigskip
{\large S. Kharchev\footnote{E-mail: kharchev@itep.ru}},\\
\bigskip
{\it Institute for Theoretical \& Experimental Physics, 117259,
Moscow,
Russia}\\
\bigskip
{\large D. Lebedev\footnote{E-mail: lebedev@mpim-bonn.mpg.de}}
\\ \bigskip
{\it Institute for Theoretical \& Experimental Physics, 117259,
Moscow, Russia} and {\it Max-Planck-Institut fur Mathematik,
Vivatsgasse 7, D-53111
Bonn, Germany},\\
\bigskip
{\large S. Oblezin} \footnote{E-mail: Sergey.Oblezin@itep.ru}\\
\bigskip {\it Institute for Theoretical \& Experimental Physics,
117259, Moscow,
Russia}\\
\end{center}

\vspace{0.5cm}

\begin{abstract}
\noindent

This paper is a short account of the construction of   a new class of the
infinite-dimensional representations of the quantum groups. The
examples include  finite-dimensional quantum
groups $U_q(\mathfrak{g})$, Yangian $Y(\mathfrak{g})$
and  affine quantum groups at zero level
 $U_q(\hat{\mathfrak{g}})_{c=0}$ corresponding to an arbitrary
 finite-dimensional semisimple Lie algebra  $\mathfrak{g}$.
At the intermediate step we construct the embedding of the quantum
groups into the algebra of the rational functions on the quantum
multi-dimensional torus.  The explicit parameterization of the quantum
groups used in this paper turns out to be closely related  to  the
parameterization of the  moduli spaces of  the monopoles. As a
result the proposed constructions of the representations provide a
quantization of the moduli spaces of  the monopoles on
$\RR^3$ and $\RR^2 \times S^1$.

\end{abstract}

\vspace{1cm}

\clearpage \newpage


\normalsize
\section{Introduction}

Constructions of  explicit realizations of  irreducible
representations are not only interesting from the general point of view
but also play a key role in the applications of the Representation
theory to Number theory, Geometry and Physics.  ``Good''
realizations often capture  important features of the underlying object
and thus lead to a  better understanding of its properties.

Below we  describe a construction of the infinite-dimensional
representations of the quantum groups that may be considered as a
generalization of the well-known construction of the
finite-dimensional representations of the classical group due to
Gelfand-Zetlin \cite{GelZ}, \cite{GelG}.  The authors came to this
construction  while trying to understand  some explicit integral
representations of the wave function of the quantum integrable
theories \cite{KL1}, \cite{KL2} obtained in the framework of Quantum Inverse
Scattering Method (QISM) \cite{F}, \cite{KS}. As it is in all other constructions
of the irreducible representations, the  crucial step is the choice of
appropriate coordinates on Lie groups and maximal commutative subalgebras in
(the skew-filed of fractions of) the universal enveloping algebra.
Late,  this construction was generalized to the
case of the quantum groups $U_q(\mathfrak{gl}(N))$ \cite{GKL3}.
The obtained realization of the representation of the
finite-dimensional quantum group possess very interesting
properties. In particular, the representation space has a natural structure
of $U_q(\mathfrak{g})\otimes U_{q^{\vee}}(\mathfrak{g}^{\vee})$~- bimodule
where $U_{q^{\vee}}(\mathfrak{g}^{\vee})$ is dual algebra. This
duality turns out to be closely related to  Langlands duality.

It would be natural to suspect that further generalizations of the
proposed constructions to the infinite-dimensional algebras are possible.
Thus in \cite{GKLO1} the construction of some class of the representations of
the Yangian $Y(\mathfrak{g})$ for an arbitrary finite-dimensional semisimple Lie
algebra $\mathfrak{g}$ was given. These representations arise as a
quantization of the symplectic leaves of the classical counterpart of the
Yangian. Surprisingly it turns out that the  symplectic leaves
of the classical Yangian  coincide with the moduli spaces of the monopoles
supplied with the appropriate symplectic structure.
The corresponding explicit expressions for the symplectic structure
on the moduli spaces of  $G$~-monopoles were  derived in \cite{AH} for $G=SU(2)$,
in \cite{B2} for $G=SU(N)$ and in \cite{FKMM} for the general case.
Thus the proposed in \cite{GKLO1} construction of the  representations of Yangian
provides at the same time the quantization of the moduli space of
monopoles. This connection between the variables
 arising in the context of QISM and the variables arising in
the study of the monopoles using the twistor methods \cite{AH}
is quite remarkable and obviously is a particular manifestation
of the deep relationship between these two subjects.

In this note we further generalize the construction of \cite{GKL1},
\cite{GKL3}, \cite{GKLO1} to obtain a realization of a class of the
infinite-dimensional representations of the finite-dimensional quantum groups
$U_q(\mathfrak{g})$ and affine quantum groups $U_q(\hat{\mathfrak{g}})_{c=0}$
with zero level $c=0$ for an arbitrary semi-simple Lie algebra $\mathfrak{g}$.
As an intermediate step, we construct the embedding of  quantum
groups into the algebra of  rational functions on the quantum
multi-dimensional torus. We postpone the technical details to another
occasion \cite{GKLO2} and concentrate on the explicit expressions for the
representations of the quantum groups.

Similar to the connection of the  Yangian representations
with the quantization of the monopoles on $\RR^3$
the proposed representations of the affine algebra
are connected with the quantization of the periodic monopoles
on $\RR^2\times S^1$. In particular the classification of the trigonometric
$r$-matrices underlying the quantum affine algebras
$U_q(\hat{\mathfrak{g}})$ \cite{BD}  corresponds to the classification of the
special class of asymptotic boundary conditions on a
monopole solutions on $\RR^2\times S^1$. It would be natural to make one
step further and consider the quantization of the moduli space of
the double-periodic monopoles on $\RR\times S^1\times
S^1$. Presumably this should correspond to the quantum elliptic
algebras and, indeed, the choice of the asymptotic boundary conditions may be
associated with elliptic $r$-matrix \cite{BD}.
The detailed account of the relevant description of the  moduli spaces of
monopoles on $\RR^2\times S^1$ and $\RR\times S^1\times  S^1$ and its
relation to the representations of the quantum groups will be given elsewhere
\cite{GKLO2}.

Finally note that the embedding of $U_q (\frak{g})$ into non-commutative
multi-dimensional torus proposed in this note differs from the known similar
constructions \cite{MV}, \cite{VK}, \cite{Nou}. Also there is an obvious
similarity of our constructions with the constructions in \cite{FO1},
\cite{FO2}. However  results presented in this note  seem to be new.

The  plan of the paper is as follows. In Section 2 we recall the
construction of a certain class of representations of the Yangian
introduced in \cite{GKLO1}. In Section 3 we give  the  generalization of
the construction discussed in Section 2 to the case of the universal
enveloping of the quantum affine algebra $U_q(\hat{\mathfrak{g}})$ at $c=0$
for an arbitrary semisimple Lie algebra $\frak{g}$. In Section 4 the explicit
construction of the representations of the quantum groups $U_q (\frak{g})$
for an arbitrary semisimple Lie algebra $\mathfrak{g}$ is given.

{\em Acknowledgments}: The  authors are grateful
 to M. Finkelberg, V. Gorbunov, A. Levin and A. Rosly for useful discussions.
The research was partly supported by grants
CRDF RM1-2545; INTAS 03-513350; grant NSh 1999.2003.2 for support of
scientific schools, and by grants RFBR-03-02-17554 (A. Gerasimov, D. Lebedev,
S. Oblezin), and RFBR-03-02-17373 (S. Kharchev).
The research of A.~Gerasimov was also partly supported by SFI Basic
Research Grant. D.~Lebedev and S. Oblezin would like to thank the Max-Planck-Institut f\"{u}r
Mathematik for  support and warm hospitality.
The research of S. Oblezin was also partly supported as a Independent
University M\"obius Prize Fellow; also S. Oblezin is thankful to the
Jumelage (Twinship) Program CNRS-IUM for partial support and to
l'Universit\'e d'Angers for the hospitality.

\section{A representation  of $Y(\frak{g}) $}
In this section we remind the explicit construction of a class of
representations of the Yangian in terms of difference operators given in
\cite{GKLO1}.

We start with  the definition of the Yangian for a semisimple
Lie algebra $\frak{g}$ due to Drinfeld  \cite{Dr1}. Let
$\mathfrak{h}\subset\mathfrak{b}\subset\mathfrak{g}$ be a simple
finite-dimensional Lie algebra $\frak{g}$ of rank $\ell$ over $\CC$ with a
fixed Cartan subalgebra $\frak{h}$ and a Borel subalgebra $\frak{b}$. Let
$a=||a_{ij}||,\,i,j=1,\ldots,\ell$ be the Cartan matrix of $\frak g$,
$\Gamma$ be the set of vertices of the Dynkin diagram of $\mathfrak{g}$,
$\{\alpha_i \in\frak{h}^{*}, i\in \Gamma\}$ be the set of simple roots and
$\{\alpha_{i}^{\vee}, i\in\Gamma\}$ be the set of the corresponding coroots
such that $a_{ij}=\alpha^{\vee}_{i}(\alpha_{j})$. There exist coprime
positive integers $d_1,\ldots,d_\ell$ such that the matrix $||d_{i}a_{ij}||$
is symmetric. Define the invariant bilinear form on $\frak{h}^{*}$ by
$(\alpha_i ,\alpha_j)=d_i a_{ij}$, then
$a_{ij}=\frac{2(\alpha_i,\alpha_j)}{(\alpha_i,\alpha_i)}$.

Introduce the formal generating  series
$H_i(u)\,,E_i(u),\,F_i(u)\,,\;i\in\Gamma$:
\be\label{int71}
H_i(u)=1+\sum_{n=0}^\infty H_i^{(n)}u^{-n-1}\;,\\
E_i(u)=\sum_{n=0}^\infty E_i^{(n)}u^{-n-1}\,,\ \ \ \ \ \
F_i(u)=\sum_{n=0}^\infty F_i^{(n)}u^{-n-1}\;.
\ee
\begin{de}
The Yangian $Y(\frak{g})$ is the associative algebra with
the elements $H^{(n)}_i, E^{(n)}_i,$ $F^{(n)}_i$, $i\in \Gamma$;
$n=0,1,2\ldots$ and the following defining relations
\be\label{cr11}
[H_i(u),H_j(v)]=0\;,\\

[H_i(u),E_j(v)]=-\frac{\imath\hbar}{2}\,(\alpha_i,\alpha_j)\,
\frac{[H_i(u),E_j(u)-E_j(v)]_{+}}{u-v}\,,\\

[H_i(u),F_j(u)]=\frac{\imath\hbar}{2}\,(\alpha_i,\alpha_j)\,
\frac{[H_i(u),F_j(u)-F_j(v)]_{+}}{u-v}\,,\\

[E_i(u),F_j(v)]=-\imath\hbar\,\frac{H_i(u)-H_i(v)}
{u-v}\,\delta_{i,j}\,,\\
\hspace{-0.5cm}
[E_i(u),E_i(v)]=
-\frac{\imath\hbar}{2}(\alpha_i,\alpha_i)\frac{(E_i(u)-E_i(v))^2}{u-v}\,,\\
\hspace{-0.5cm}
[F_i(u),F_i(v)]=
\frac{\imath\hbar}{2}(\alpha_i,\alpha_i)\frac{(F_i(u)-F_i(v))^2}{u-v}\,,\\
\begin{array}{l}
\hspace{4cm}[E_i(u),E_j(v)]=\\-\frac{\imath\hbar}{2}
(\alpha_i,\alpha_j)\frac{[E_i(u),E_j(u)-E_j(v)]_+}{u-v}-
\frac{[E_i^{(0)},E_j(u)-E_j(v)]}{u-v}\,,\\
\hspace{4cm}
[F_i(u),F_j(v)]=\\ \frac{\imath\hbar}{2}
(\alpha_i,\alpha_j)\frac{[F_i(u),F_j(u)-F_j(v)]_+}{u-v}-
\frac{[F_i^{(0)},F_j(u)-F_j(v)]}{u-v}\,,\\
\end{array}\ \ \ \ i\neq j\,,
\ee
\be\label{cr51}
\sum_{\sigma\in\frak S_m} [E_i(u_{\sigma(1)}),[E_i(u_{\sigma(2)}),
\ldots,[E_i(u_{\sigma(m)}),E_j(v)]\ldots]]=0\;,\\
\sum_{\sigma\in\frak S_m}[F_i(u_{\sigma(1)}),[F_i(u_{\sigma(2)}),
\ldots,[F_i(u_{\sigma(m)}),F_j(u)]\ldots]]=0\;,\\
m=1-a_{ij}\ \ \ \mbox{for}\ \ i\neq j\,,
\ee
where  $[a,b]_+:= ab+ba$ and summation in (\ref{cr51}) is performed over
permutation group $\frak S_m$.
\end{de}
Let  $Y(\frak{b})\subset Y(\frak{g})$ be the subalgebra generated
by $H_i(u), E_i(u),\,i\in\Gamma$.

The explicit description of the representation of the Yangian in terms of
difference ope\-rators is based  on the choice of a large enough commutative
subalgebra. We shall use the coefficients of the series $H_i(u)$ as the
generators of this subalgebra. In the constructed representation $H_i^{(n)}$
will act by multiplication on the  functions of some auxiliary variables.
Thus to obtain representation of the Yangian we should find the representation
of the other generators in terms of some  difference operators acting on the
same space  of  functions.

Let us introduce a set of variables
$\{\gamma_{i,k}\,\,;i\in\Gamma; \,k=1,\ldots,m_i\}$, where $m_i\in\NN$
and let $\cal M$ be the space of meromorphic functions in these variables.
Define  the following difference operators acting on $\cal M$:
$\beta_{i,k}=e^{\i\hbar d_i\frac{\d}{\d\gamma_{i,k}}}$.
Below we use the convention $\prod_{s=j}^k f_s:=1$, for any $f_s$ if $k<j$.

Consider the operators
\be\label{mn5}
H_i(u)=R_i(u)\frac{\prod\limits_{j\neq i}\prod\limits_{r=1}^{-a_{ji}}
\prod\limits_{p=1}^{m_j}(u-\gamma_{j,p}-\ts{\i\h}{2}(\al_i+r\al_j,\al_j))}
{\prod\limits_{p=1}^{m_i}(u-\gamma_{i,p})(u-\gamma_{i,p}-\ts{\i\h}{2}
(\al_i,\al_i))}\,,
\ee
\be\label{mn6}
E_i(u)=d_i^{-1/2}\sum\limits_{k=1}^{m_i}
\frac{\prod\limits_{j=i+1}^\ell\prod\limits_{r=1}^{-a_{ji}}
\prod\limits_{p=1}^{m_j}(\gamma_{i,k}-\gamma_{j,p}-
\ts{\i\h}{2}(\al_i+r\al_j,\al_j))}
{(u-\gamma_{i,k})\prod\limits_{p\neq k}(\gamma_{i,k}-\gamma_{i,p})}\,
\beta_{i,k}^{-1}\,,
\ee
\be\label{mn7}
F_i(u)=-d_i^{-1/2}\sum\limits_{k=1}^{m_i}R_i(\gamma_{i,k}+
\ts{\imath\hbar}{2}(\alpha_i,\alpha_i))\;\times\\
\frac{\prod\limits_{j=1}^{i-1}\prod\limits_{r=1}^{-a_{ji}}
\prod\limits_{p=1}^{m_j}(\gamma_{i,k}-\gamma_{j,p}-
\ts{\i\h}{2}(\al_i+r\al_j,\al_j)+\ts{\i\h}{2}(\al_i,\al_i))}
{(u-\gamma_{i,k}-\frac{\imath\hbar}{2}
(\alpha_i,\alpha_i))\prod\limits_{p\neq k}(\gamma_{i,k}-\gamma_{i,p})}
\,\beta_{i,k}\,,
\ee
where the rational functions $R_i(u)$ will be specified below.

\begin{te}\label{MainTH}
\noindent \cite{GKLO1} (i). For any set of positive integers $\{m_i,\,,
i\in\Gamma\}$
satisfying the conditions $l_i:=\sum_{j=1}^{\ell}m_j a_{ji}\in\ZZ_{\geq 0}$,
introduce the polynomials $R_i(u)=\prod_{s=1}^{l_i}(u-\nu_{i,s})$,
where $\{\nu_{i,s}\,,i\in\Gamma$, $s=1,\ldots,l_i\}$
is a set of arbitrary complex parameters.
Then the operators (\ref{mn5})-(\ref{mn7}) considered as formal power
series in $u^{-1}$, define a representation of $\,Y(\frak{g})$ in the space
$\cal M$.

\noindent (ii). Let $\{m_i,\,i\in\Gamma\}$ be an  arbitrary set of
positive integers and $R_i(u)$ be rational functions of the form
$R_i(u)=\prod_{s=1}^{l^+_i}(u-\nu^+_{i,s})\big/
\prod_{s=1}^{l^-_i}(u-\nu^-_{i,s})$,
where $\{\nu^{\pm}_{i,s}\,,i\in\Gamma\,,s=1,\ldots,l^{\pm}_i\}$
is a set of arbitrary complex parameters and
$\,l_i^+-l_i^-=\sum_{j=1}^{\ell}m_ja_{ji}$. Then the operators
(\ref{mn5}), (\ref{mn6}) considered as formal power
series in $u^{-1}$, define a representation of $\,Y(\frak{b})$ in the space
$\cal M$.
\end{te}

\noindent
Below we generalize this construction to the case of the quantum
affine algebra at zero level.

\section{A representation of $U_q(\hat{\mathfrak{g}})_{c=0}$}

In this section we construct a representation of
 the quantum affine algebra at the zero level  $U_q({\hat{\frak{g}}})_{c=0}$.
 The construction is the direct
generalization of the  construction for the Yangian described  in the previous section.
We start with  the definition of the quantum affine algebra
$U_q(\hat{\frak{g}})$ at $c=0$  for any semisimple
Lie algebra $\frak{g}$ in terms of generating series following
 \cite{Dr1}.

Let $q$ be an undeterminate.
Quantum affine algebra as an  associative $\CC(q)$-algebra may be
present in terms of the elements $K_i^{\pm 1},\,H_i^{(n)},\,
n\in\ZZ\backslash\{0\},\, E_i^{(n)},F_i^{(n)}, n\in\ZZ$, $i\in\Gamma$ and
relations. Introduce the formal generating series $K_i^{\pm} (z)\,,
E_i (z),\,$ and $ F_i(z)\,,\;i\in\Gamma$:
\be\label{int7}
K_i^{\pm}(z)=K_i^{\pm 1}\exp\big(\pm(q_i-q_i^{-1})\sum_{n\in\NN}H_i^{(\pm n)}
z^{\mp n}\big)\;,\\
E_i(z)=\sum_{n\in\ZZ}E_i^{(n)}z^{-n}\,,\ \ \ \ \ \
F_i(z)=\sum_{n\in\ZZ}F_i^{(n)}z^{-n}\;,
\ee
where $q_i:=q^{d_i}$.
\begin{de} Quantum affine algebra
  $U_q (\hat{\frak{g}})$ at ${c=0}$ is the associative algebra with
elements $K_i^{\pm 1}, H^{(n)}_i$, $n\in\ZZ\backslash\{0\}$;
$E^{(n)}_i, F^{(n)}_i$, $n\in\ZZ$; $i\in\Gamma$ and the following defining
relations:
\be\label{cr1}
K_i^{\pm}(z)K_j^{\pm}(w)=K_j^{\pm}(w)K_i^{\pm}(z),\\
K_i^{+}(z)K_j^{-}(w)=K_j^{-}(w)K_i^{+}(z)\,,\\
(z-q_i^{a_{ij}}w)K_i^{\pm}(z)E_j(w)=(q_i^{a_{ij}}z-w)
E_j(w)K^{\pm}_i(z)\,,\\
(z-q_i^{-a_{ij}}w)K_i^{\pm}(z)F_j(w)=(q_i^{-a_{ij}}z-w)F_j(w)K^{\pm}_i(z)\,,\\

[E_i(z),F_j(w)]=\frac{\delta_{i,j}}{q_i-q_i^{-1}}\,\delta(z/w)\big(K^{+}_i
(w)-K^{-}_i(z)\big),\\
(z-q_i^{a_{ij}}w)E_i(z)E_j(w)=(q_i^{a_{ij}}z-w)E_{j}(w)E_{i}(z),\\
(z-q_i^{-a_{ij}}w)F_i(z)F_j(w)=(q_i^{-a_{ij}}z-w)F_{j}(w)F_i(z)\,,\\
\hspace{-0.5cm}
\sum_{\sigma\in{\frak S}_m}\sum_{k=0}^m(-1)^k
\genfrac{[}{]}{0pt}{0}{m}{k}_{q_i}E_i(z_{\sigma(1)})
\ldots E_i (z_{\sigma(k)})E_j(w)E_i(z_{\sigma(k+1)})\ldots
E_i(z_{\sigma(m)})=0\;,\hspace{-0.2cm}\\
\hspace{-0.5cm}
\sum_{\sigma\in{\frak Sm}}\sum_{k=0}^m(-1)^k
\genfrac{[}{]}{0pt}{0}{m}{k}_{q_i}F_i(z_{\sigma(1)})
\ldots F_i (z_{\sigma(k)})E_j(w)F_i(z_{\sigma(k+1)})\ldots
F_i(z_{\sigma(m)})=0\;,\\
m=1-a_{ij}\ \ \ \mbox{for}\ \ i\neq j\,.
\ee
\end{de}
Here we use the standard notations
$\displaystyle
\genfrac{[}{]}{0pt}{0}{m}{k}_q=\frac{[m]_{q}!}{[k]_{q}![m-k]_{q}!}\,,\;
[k]_{q}!=\prod\limits_{1\leq j \leq k}\frac{q^j-q^{-j}}{q-q^{-1}}\,.
$
The formal delta-function is defined as $\delta(z)=\sum_{n\in\ZZ}z^n$.
For more details on the operator-valued formal series see, for example,
\cite{K}.

To define a representation of $U_q(\hat{\frak g})_{c=9}$ in terms of the
difference operators, we start with construction of the embedding of the
corresponding universal enveloping algebra into the algebra $\TT_q$
of the rational functions of the non-commutative multi-dimensional torus.

Let $\TT_q$ be the associative $\CC(q)$-algebra
of the rational functions of invertible elements
${\bf v}_{i,k},\,{\bf u}_{i,k},\,{\bf w}_{i,s},\,
i=1,\ldots,\ell;\,k=1,\ldots,m_i\in\NN;\,s=1,\ldots,l_i\in\ZZ_{\geq0}$,
subject to relations
\be\label{tq1}
{\bf v}_{i,k}{\bf v}_{j,l}={\bf v}_{j,l}{\bf v}_{i,k}\,,\ \ \ \
{\bf u}_{i,k}{\bf u}_{j,l}={\bf u}_{j,l}{\bf u}_{i,k}\,,\\
{\bf u}_{i,k}{\bf v}_{j,l}=q_i^{\delta_{i,j}\delta_{k,l}}{\bf v}_{j,l}
{\bf u}_{i,k}\,,
\ee
and ${\bf w}_{i,s}$ are central elements in $\TT_q$.
Chose the set of natural numbers $\{ m_i\}$ satisfying the conditions
\be\label{cond}
\sum\limits_{j=1}^{\ell}m_j a_{ji}=l_i,
\ee
and consider the following formal generating functions  $K_i^{\pm}$,
$E_i(z)$ and $F_i(z)$ in variable $z$.

The functions $K^+_i (z)$ and $K^-_i (z)$ are defined as  infinite series
expansion in $z^{-1}$ and $z$   of the same rational function
\be\label{grep1}
K_i(z)=c_i
\prod\limits_{j=1}^{\ell}\prod\limits_{p=1}^{m_j}{\bf v}_{j,p}^{a_{ji}}
\cdot\prod_{s=1}^{l_i}(z{\bf w}_{i,s}^{-1}-{\bf w}_{i,s})\;
\frac{\prod\limits_{j\neq i}\prod\limits_{r=1}^{-a_{ji}}
\prod\limits_{p=1}^{m_j}(z-q_j^{a_{ji}+2r}{\bf v}_{j,p}^2)}
{\prod\limits_{p=1}^{m_i}(z-{\bf v}_{i,p}^2)(z-q_i^{2}{\bf v}_{i,p}^2)}\,,
\ee
where $c_i=\prod_{j=1}^\ell q_j^{m_ja_{ji}/2}$.
Due to conditions (\ref{cond}) the degrees of the polynomials in the
numerator and  in the denominator in  (\ref{grep1}) coincide and thus
$K_i^{\pm}(z)$ as the power series in $z^{\mp 1}$  have the form (\ref{int7})
where $K_i=c_i\prod_{s=1}^{l_i}{\bf w}_{i,s}^{-1}\cdot
\prod_{j=1}^{\ell}\prod_{p=1}^{m_j}{\bf v}_{j,p}^{a_{ji}}$.

Let $R^{(\pm)}(z)$ be the polynomials such that
$R^{(+)}_i(z)R^{(-)}_i(z)=\prod_{s=1}^{l_i}(z{\bf w}_{i,s}^{-1}-
{\bf w}_{i,s})$. The other generating series are defined as follows
\be\label{grep2}
E_i(z)=\frac{c_i}{q_i-q_i^{-1}}
\prod\limits_{p=1}^{m_i}{\bf v}_{i,p}\cdot\prod\limits_{j=i+1}^{\ell}
\prod\limits_{p=1}^{m_j}{\bf v}_{j,p}^{a_{ji}}\;\cdot\\
\sum\limits_{k=1}^{m_i}{\delta}(z/{\bf v}_{i,k}^2)\,
{\bf v}_{i,k}^{-2}R_i^{(+)}({\bf v}_{ik}^2)
\frac{\prod\limits_{j=i+1}^{\ell}\prod\limits_{r=1}^{-a_{ji}}
\prod\limits_{p=1}^{m_j}({\bf v}_{i,k}^{2}-q_j^{a_{ji}+2r}
{\bf v}_{j,p}^2)}{\prod\limits_{p\neq k}
({\bf v}_{i,k}^{2}-{\bf v}_{i,p}^2)}\;{\bf u}_{i,k}^{-1}\,,\\
F_i(z)=-\frac{q_i^{-2m_i}}{q_i-q_i^{-1}}
\prod\limits_{p=1}^{m_i}{\bf v}_{i,p}\cdot
\prod\limits_{j=1}^{i-1}\prod\limits_{p=1}^{m_j}{\bf v}_{j,p}^{a_{ji}}\;
\cdot\\
\sum\limits_{k=1}^{m_i}{\delta}(z/q_i^{2}{\bf v}_{i,k}^2)\,
{\bf v}_{i,k}^{-2}R_i^{(-)}(q_i^{2}{\bf v}_{i,k}^2)
\frac{\prod\limits_{j=1}^{i-1}\prod\limits_{r=1}^{-a_{ji}}\prod
\limits_{p=1}^{m_j}(q_i^2{\bf v}_{i,k}^{2}-q_j^{a_{ji}+2r}{\bf v}_{j,p}^2)}
{\prod\limits_{p\neq k}({\bf v}_{i,k}^{2}-{\bf v}_{i,p}^2)}\;{\bf u}_{i,k}\,.
\ee
\begin{te}\label{MTH}
For any set of positive integers $\{ m_i, i\in\Gamma\}$ obeying the
conditions
(\ref{cond}), the generating series $K^{\pm}_i(z)$, $E_i(z)$, and
$F_i(z)$ defined by (\ref{grep1}), (\ref{grep2}), satisfy the relations
(\ref{cr1}) and therefore define an embedding $\pi: U_q(\hat
{\frak{g}})_{c=0}$ $\hookrightarrow\TT_q$.
\end{te}
The proof will be given in \cite{GKLO2}.

Let us introduce a set of variables
$\{\gamma_{i,k}\in\CC\,;i\in\Gamma; \,k=1,\ldots,m_i\}$, and let $\cal M$
be the space of meromorphic functions in
these variables. We also fix a set of the complex numbers $\nu_{i,s}$,
$i=1,\ldots,\ell,\,s=1,\ldots,l_i$.
Define the following representation of $\TT_q$ in terms of
difference operators acting on $\cal
M$:
\be\label{shift}
{\bf u}_{i,k}=e^{\i\omega_1 d_i\frac{\d}{\d\gamma_{i,k}}},\,\,\,\,
{\bf v}_{i,k}=e^{\frac{2\pi\gamma_{ik}}{\omega_2}},\,\,\,\,
{\bf w}_{i,s}=e^{\frac{2\pi\nu_{i,s}}{\omega_2}},\,\,\,\,
q=e^{\frac{2\pi\i\omega_1}{\omega_2}},
\ee
where $\omega_1$ and $\omega_2$ are arbitrary complex parameters.
The following proposition is the simple consequence of Theorem \ref{MTH}.

\begin{prop} Let us given the representation of $~\TT_q$ defined by
(\ref{shift}). The coefficients $K_i^{\pm 1},\,H_i^{(n)},\,
n\in\ZZ\backslash\{0\},\, E_i^{(n)},F_i^{(n)}, n\in\ZZ$, $i\in\Gamma$ of
the formal power series (\ref{grep1}), (\ref{grep2}) define a representation
of $U_q(\hat{\frak{g}})_{c=0}$ in $\cal M$.
\end{prop}

\section{A representation of $U_q(\frak{g})$}

The finite-dimensional quantum groups  $U_q(\frak{g})$ may be naturally
considered as the subalgebras of the affine quantum groups
$U_q(\hat{\frak{g}})_{c=0}$.
Indeed, the generators $K_i^{\pm 1}$, $E^{(0)}_i,$ and $F^{(0)}_i$ introduced
in the previous section  obey the commutation relations of the
finite-dimensional quantum group $U_q(\frak{g})$. We shall omit the
superscript $(0)$ and write down simply $E_i$ and $F_i$.
Thus the representations of the
affine quantum groups introduced above automatically provide the
representations of the finite-dimensional quantum groups. Below we give the
explicit expressions for the generators in this representations. However we
would prefer to work in this section in a slightly more general case
and consider the various rational forms of the quantum groups.
Let $Q$ and $P$ be  the root and weight lattices of the Lie algebra
$\mathfrak{g}$.  There are different rational forms of the quantum
group not isomorphic as $\CC(q)$ algebras which may be enumerated
by the choice of a sublattice $M$ such that $Q\subseteq M\subseteq P$
(see \cite{CP}, \cite{L}
and references therein). Thus for example the adjoint (the smallest)
rational form $U_q^{Q}(\frak{g})$ is the associative $\CC(q)$-algebra
generated by $E_i ,F_i$ and $K^{\pm}_i$, $i\in\Gamma$
with the defining relations:

\vspace{-0.3cm}
\be\label{def1}
K_{i}K_{i}^{-1}=K_{i}^{-1}K_{i}=1\,,\ \ \ \ K_{i}K_{j}=K_{j}K_{i}\,,\\

\vspace{-0.1cm}
K_{i}E_{j}K_{i}^{-1}=q_i^{a_{ij}}E_{j}\,,\ \ \ \
K_{i}F_{j}K_{i}^{-1}=q_i^{-a_{ij}}F_{j}\,,\\

\vspace{-0.2cm}
E_{i}F_{j}-F_{j}E_{i}=\delta_{i,j}
\frac{K_i -K^{-1}_i}{q_i-q_i^{-1}}\,,\\

\vspace{-0.2cm}
\sum\limits_{r=0}^{1-a_{ij}}(-1)^{r}\genfrac{[}{]}{0pt}{0}{1-a_{ij}}{r}_{q_i}
E_{i}^{1-a_{ij}-r} E_{j}E_{i}^{r}=0\,,\ \ \ \ \ i\neq j\,,\\
\sum\limits_{r=0}^{1-a_{ij}}(-1)^{r}\genfrac{[}{]}{0pt}{0}{1-a_{ij}}{r}_{q_i}
F_{i}^{1-a_{ij}-r} F_{j}F_{i}^{r}=0\,, \ \ \ \ \ i\neq j\,.
\ee
On the other hand, the largest, simply-connected, rational from
$U^{P}_{q}(\frak{g})$ is obtained by adjoining to
$U_{q}^{Q}(\frak{g})$ the invertible elements $L_i
\,,i\in\Gamma,$ such that $K_i =\prod_j L_{j}^{a_{ji}}$ and obeying the
relations $L_{i}E_{j}L_{i}^{-1}=q_i^{\delta_{i,j}}E_{j},\,
L_{i}F_{j}L_{i}^{-1}=q_i^{-\delta_{i,j}}F_{j}$.
Given any intermediate lattice  $Q \subseteq M\subseteq P$, the
corresponding rational form  $U_{q}^{M}({\frak g})$ is obtained by
adjoining to $U^{Q}_{q}({\frak g})$ the elements
$K_{\beta_i}=\prod_{j}L_{j}^{n_{ji}} $ for any basis
$\beta_i =\sum_{j}n_{ji}\lambda_{j}\in M,\,i\in\Gamma$,
where $\lambda_j$ are the fundamental weights. For any simple root
$\alpha_i$  let $\alpha_i=\sum_{j}m_{ji}\beta_{j}$. Then the generators
$K_i$ are expressed by $K_{\b_j}$ as  $K_{i}=\prod_{j}K_{\beta_j}^{m_{ji}}$.
Let $||M_{ij}||$ be the matrix
inverse to $||m_{ij}||$ and $d=\det||a_{ij}||$. Let $\mathbb{T}_q$
be the quantum torus defined in the  previous
section. The following Theorem is a generalization of the
embedding of the finite-dimen\-si\-onal  quantum group into quantum torus
obtained by the restriction of embedding  of the affine quantum group
described in Section 3.
\begin{te}
(i). Let $\frak{g}$ be an arbitrary semisimple Lie algebra.
Let $K_{\beta_i}^{\pm 1},E_i,F_i$, $i\in\Gamma$
be the generators of the rational form of quantum group $U^{M}_q(\frak{g})$,
associated with the lattice $M$ such that $Q\subseteq M\subseteq P$.
The following  expressions define the embedding
$\pi: U^{M}_q(\frak{g})\hookrightarrow\TT_q:$
\be\label{rg1}
\pi(K_{\beta_i})=c_{\beta_i}\prod\limits_{j=1}^{\ell}\Big(
\prod_{s=1}^{l_j}{\bf w}_{j,s}^{-dM_{ji}}\cdot\prod\limits_{p=1}^{m_j}
{\bf v}_{j,p}^{n_{ji}}\Big)\,,\\
\pi(E_i)=\frac{c_{\al_i}}{q_i-q_i^{-1}}
\prod\limits_{p=1}^{m_i}{\bf v}_{i,p}\cdot\prod\limits_{j=i+1}^{\ell}
\prod\limits_{p=1}^{m_j}{\bf v}_{j,p}^{a_{ji}}\;\cdot\\
\sum_{k=1}^{m_i}{\bf v}_{i,k}^{-2}R_i^{(+)}({\bf v}_{ik}^2)
\frac{\prod\limits_{j=i+1}^{\ell}\prod\limits_{r=1}^{-a_{ji}}
\prod\limits_{p=1}^{m_j}({\bf v}_{i,k}^{2}-q_j^{a_{ji}+2r}
{\bf v}_{j,p}^2)}{\prod\limits_{p\neq k}
({\bf v}_{i,k}^{2}-{\bf v}_{i,p}^2)}\;{\bf u}_{i,k}^{-1}\,,\\
\pi(F_i)=-\frac{q_i^{-2m_i}}{q_i-q_i^{-1}}
\prod\limits_{p=1}^{m_i}{\bf v}_{i,p}\cdot
\prod\limits_{j=1}^{i-1}\prod\limits_{p=1}^{m_j}{\bf v}_{j,p}^{a_{ji}}\;
\cdot\\
\sum_{k=1}^{m_i}{\bf v}_{i,k}^{-2}R_i^{(-)}(q_i^{2}{\bf v}_{i,k}^2)
\frac{\prod\limits_{j=1}^{i-1}\prod\limits_{r=1}^{-a_{ji}}\prod
\limits_{p=1}^{m_j}(q_i^2{\bf v}_{i,k}^{2}-
q_j^{a_{ji}+2r}{\bf v}_{j,p}^2)}
{\prod\limits_{p\neq k}({\bf v}_{i,k}^{2}-{\bf v}_{i,p}^2)}\;{\bf u}_{i,k}\,,
\ee
where $R^{(\pm)}_i(z)$ be the polynomials such that
$R_i^{(+)}(z)R_i^{(-)}(z)=\prod_{s=1}^{l_i}
(z{\bf w}_{i,s}^{-d}-{\bf w}_{i,s}^{d})$ and
$c_{\beta_i}=\prod_{j=1}^\ell q_j^{m_jn_{ji}/2}$.

\noindent
(ii). Let the elements of $~\TT_q$ are represented by (\ref{shift}).
Then the operators (\ref{rg1}) define a representation of $U^M_q({\frak{g}})$
in ${\cal M}$.
\end{te}
\begin{rem}
In the case $\frak{g}=\frak{sl}(\ell+1)$  the Theorem 4.1
is a natural generalization of  the  Theorem 3.1 in \cite{GKL3}.
\end{rem}

Note that there are known several embeddings of $U_q(\frak{g})$ into the
algebra of quantum torus ( e.g. \cite{Nou}, \cite{MV}, \cite{VK}). Our
construction is the new one and has deep relations with the Quantum Inverse
Scattering Method \cite{GKL1} as well as with the natural parameterization of
the moduli spaces of monopoles \cite{GKLO1}.
The detailed discussion of the connection with the monopoles on
$\RR^2\times S^1$ will be described in \cite{GKLO2}.


\begin{thebibliography}{100}

\bibitem{AH}
M.F. Atiyah,  N. Hitchin,  {\it The geometry and dynamics of
magnetic monopoles},  Prinston, NJ University Press (1988).

\bibitem{Nou}
H. Awata, M. Noumi, S. Odake, {\it Heisenberg realization for
$U_q(sl_n)$ on the flag manifold}, Lett. Math. Phys. {\bf 30}
No.1, (1994), 35--43.

\bibitem{BD}
A. Belavin and V. Drinfeld {\it  Triangle equations and simple Lie algebras},
Funct. Anal. Appl., {\bf 16} (1982), 159.

\bibitem{B2}
R. Bielawski {\it Asymptotic metrics for $SU(N)$- monopoles with maximal
symmetry breaking}. Comm. Math. Phys. {\bf 199} (1998), 297--325.

\bibitem{CP}
V. Chari, A. Pressley, {\it A guide to quantum groups},
Cambridge Univ. Press, Cambridge, 1994.

\bibitem{Dr0}
V.G. Drinfeld, {\it Hopf algebras and the quantum Yang-Baxter equation},
Soviet Math. Dokl. {\bf 32}, 254--258, (1985).

\bibitem{Dr1}
V.G. Drinfeld, {\it A new realization of Yangians and of
quantum affine algebras}, Dokl. Akad. Nauk SSSR {\bf 296} (1987), no. 1,
13--17 (Russian); translated in Soviet Math. Dokl. {\bf 36} (1988), 212--216.

\bibitem{F}
L.D. Faddeev, {\it Quantum completely integrable models in field
theory}, Sov. Sci. Rev., Sect. C (Math. Phys. Rev.) {\bf 1} (1980), 107--155.

\bibitem{FO1}
B. Feigin, A. Odesskii, {\it Vector bundles on
Elliptic Curves and Sklyanin Algebras} in: {\it Topics in quantum
groups and finite type invariants, Mathematics at the Independent
University of Moskow, B. Feigin and V. Vasiliev eds,Advances in Mathematical
Sciences} {\bf 38}, AMS Translations, ser. 2, {\bf 185} (1998), 65--84.

\bibitem{FO2}
B. Feigin, A. Odesskii, {\it Elliptic deformations of current algebras and
their representations by difference operators}, Funk. Anal. Appl. {\bf 31}
(1997), 57--70.

\bibitem{FKMM}
M. Finkelberg, A.  Kuznetsov, N.  Markarian, I. Mirkovi\'{c}, {\it  A note on
the symplectic structure on the space of G-monopoles},
Comm. Math. Phys. {\bf 201} (1999), 411--421.

\bibitem{GelZ}
I.M. Gelfand, M.L. Tsetlin, {\it Finite-dimensional representations of the
group of uni\-modular matrices}, Dokl.Akad. Nauk SSSR {\bf 71} (1950),
825--828 (Russian), translated in I.M. Gelfand,
{\it Collected Papers, Vol. II}, Springer, Berlin, 1988, 653--656.

\bibitem{GelG}
I.M. Gelfand, M.I. Graev, {\it Finite-dimensional
irreducible representations of the unitary and the full linear groups,
and related special functions}, Izv. Akad. Nauk SSSR, Ser.Mat. {\bf 29}
(1965), 1329--1356; translated in Amer. Math. Soc. Trans.
Ser. 2 {\bf 64} (1965), 116--146.

\bibitem{GKL1}
A. Gerasimov, S. Kharchev, D. Lebedev,
{\it Representation Theory and Quantum Inverse Scattering Method:
The Open Toda Chain and the Hyperbolic Sutherland Model,}
Int. Math. Res. Notices {\bf 17} (2004), 823--854.

\bibitem{GKL2}
A. Gerasimov, S. Kharchev, D. Lebedev, {\it  On a class of integrable
systems connected with $GL(N,\RR)$,} Int. J. Mod. Phys. A {\bf 19} Suppl.,
(2004), 205--216.

\bibitem{GKL3}
A. Gerasimov, S. Kharchev,  D. Lebedev, {\it Representation theory and
quantum integrability}, arXiv: math.QA/0402112.

\bibitem{GKLO1}
A. Gerasimov, S. Kharchev, D. Lebedev,  S. Oblezin, {\it On a class of
representations of the Yangian and moduli space of  monopoles},
arXiv: math.AG/0409031.

\bibitem{GKLO2}
A. Gerasimov, S. Kharchev, D. Lebedev,  S. Oblezin, in preparation.



\bibitem{K}
V. Kac, {\it Infinite dimensional Lie algebras}, Third edition,
Cambridge University Press, Cambridge, 1990.

\bibitem{KL1}
S. Kharchev, D. Lebedev, {\it Eigenfunctions of $GL(N,\RR)$ Toda chain: The
Mellin-Barnes representation}, JETP Lett. {\bf 71} (2000), 235--238.

\bibitem{KL2}
S. Kharchev, D. Lebedev, {\it Integral representations for the eigenfunctions
of quantum open and periodic Toda chains from QISM formalism},
J. Phys. {\bf A34} (2001), 2247--2258.

\bibitem{KS}
P.P. Kulish, E.K. Sklyanin, {\it Quantum spectral transform
method. Recent developments}, Lecture Notes in Phys.
{\bf 151}, pp. 61--119, Springer, Berlin-New York, 1982.

\bibitem{L}
G. Lusztig, {\it Introduction to quantum groups},
Progress in Mathematics, 110, Birkh\"auser Boston, Inc., Boston, MA, 1993.

\bibitem{MV}
A. Morozov and L.  Vinet, {\it Free-Field Representation  of Group Element
for Simple Quantum Group},  Int. J. Mod. Phys. {\bf A13} (1998), 1651--1708.

\bibitem{VK}
A. Volkov and A. Kashaev {\it From the Tetrahedron Equation to Universal
R-Matrices},  L. D. Faddeev's Seminar on Mathematical Physics, 79--89,
Amer. Math. Soc. Transl. Ser. 2, 201, Amer.  Math. Soc., Providence, RI, 2000.

\end{thebibliography}
\end{document}